\magnification=1200
%\nopagenumbers

\hsize=11.25cm    
\vsize=18cm       
\parindent=12pt   \parskip=5pt     

\hoffset=.5cm   
\voffset=.8cm   

\pretolerance=500 \tolerance=1000  \brokenpenalty=5000

\catcode`\@=11

\font\eightrm=cmr8         \font\eighti=cmmi8
\font\eightsy=cmsy8        \font\eightbf=cmbx8
\font\eighttt=cmtt8        \font\eightit=cmti8
\font\eightsl=cmsl8        \font\sixrm=cmr6
\font\sixi=cmmi6           \font\sixsy=cmsy6
\font\sixbf=cmbx6

\font\tengoth=eufm10 
\font\eightgoth=eufm8  
\font\sevengoth=eufm7      
\font\sixgoth=eufm6        \font\fivegoth=eufm5

\skewchar\eighti='177 \skewchar\sixi='177
\skewchar\eightsy='60 \skewchar\sixsy='60

\newfam\gothfam           \newfam\bboardfam

\def\tenpoint{
  \textfont0=\tenrm \scriptfont0=\sevenrm \scriptscriptfont0=\fiverm
  \def\rm{\fam\z@\tenrm}
  \textfont1=\teni  \scriptfont1=\seveni  \scriptscriptfont1=\fivei
  \def\oldstyle{\fam\@ne\teni}\let\old=\oldstyle
  \textfont2=\tensy \scriptfont2=\sevensy \scriptscriptfont2=\fivesy
  \textfont\gothfam=\tengoth \scriptfont\gothfam=\sevengoth
  \scriptscriptfont\gothfam=\fivegoth
  \def\goth{\fam\gothfam\tengoth}
  
  \textfont\itfam=\tenit
  \def\it{\fam\itfam\tenit}
  \textfont\slfam=\tensl
  \def\sl{\fam\slfam\tensl}
  \textfont\bffam=\tenbf \scriptfont\bffam=\sevenbf
  \scriptscriptfont\bffam=\fivebf
  \def\bf{\fam\bffam\tenbf}
  \textfont\ttfam=\tentt
  \def\tt{\fam\ttfam\tentt}
  \abovedisplayskip=12pt plus 3pt minus 9pt
  \belowdisplayskip=\abovedisplayskip
  \abovedisplayshortskip=0pt plus 3pt
  \belowdisplayshortskip=4pt plus 3pt 
  \smallskipamount=3pt plus 1pt minus 1pt
  \medskipamount=6pt plus 2pt minus 2pt
  \bigskipamount=12pt plus 4pt minus 4pt
  \normalbaselineskip=12pt
  \setbox\strutbox=\hbox{\vrule height8.5pt depth3.5pt width0pt}
  \let\bigf@nt=\tenrm       \let\smallf@nt=\sevenrm
  \normalbaselines\rm}

\def\eightpoint{
  \textfont0=\eightrm \scriptfont0=\sixrm \scriptscriptfont0=\fiverm
  \def\rm{\fam\z@\eightrm}
  \textfont1=\eighti  \scriptfont1=\sixi  \scriptscriptfont1=\fivei
  \def\oldstyle{\fam\@ne\eighti}\let\old=\oldstyle
  \textfont2=\eightsy \scriptfont2=\sixsy \scriptscriptfont2=\fivesy
  \textfont\gothfam=\eightgoth \scriptfont\gothfam=\sixgoth
  \scriptscriptfont\gothfam=\fivegoth
  \def\goth{\fam\gothfam\eightgoth}
  
  \textfont\itfam=\eightit
  \def\it{\fam\itfam\eightit}
  \textfont\slfam=\eightsl
  \def\sl{\fam\slfam\eightsl}
  \textfont\bffam=\eightbf \scriptfont\bffam=\sixbf
  \scriptscriptfont\bffam=\fivebf
  \def\bf{\fam\bffam\eightbf}
  \textfont\ttfam=\eighttt
  \def\tt{\fam\ttfam\eighttt}
  \abovedisplayskip=9pt plus 3pt minus 9pt
  \belowdisplayskip=\abovedisplayskip
  \abovedisplayshortskip=0pt plus 3pt
  \belowdisplayshortskip=3pt plus 3pt 
  \smallskipamount=2pt plus 1pt minus 1pt
  \medskipamount=4pt plus 2pt minus 1pt
  \bigskipamount=9pt plus 3pt minus 3pt
  \normalbaselineskip=9pt
  \setbox\strutbox=\hbox{\vrule height7pt depth2pt width0pt}
  \let\bigf@nt=\eightrm     \let\smallf@nt=\sixrm
  \normalbaselines\rm}

\tenpoint

\def\pc#1{\bigf@nt#1\smallf@nt}         \def\pd#1 {{\pc#1} }

\catcode`\;=\active
\def;{\relax\ifhmode\ifdim\lastskip>\z@\unskip\fi
\kern\fontdimen2  -1.2 \fontdimen3 \string;}

\catcode`\:=\active
\def:{\relax\ifhmode\ifdim\lastskip>\z@\unskip\fi\penalty\@M\ \fi\string:}

\catcode`\!=\active
\def!{\relax\ifhmode\ifdim\lastskip>\z@
\unskip\fi\kern\fontdimen2  -1.1 \fontdimen3 \string!}

\catcode`\?=\active
\def?{\relax\ifhmode\ifdim\lastskip>\z@
\unskip\fi\kern\fontdimen2  -1.1 \fontdimen3 \string?}

\frenchspacing

\def\raggedbottom{\topskip 10pt plus 36pt\r@ggedbottomtrue}

\def\pointir{\unskip . --- \ignorespaces}

\def\Medbreak{\vskip-\lastskip\medbreak}

\long\def\th#1 #2\enonce#3\endth{
   \Medbreak\noindent
   {\pc#1} {#2\unskip}\pointir{\it #3}\smallskip}

\def\decale#1{\smallbreak\hskip 28pt\llap{#1}\kern 5pt}
\def\decaledecale#1{\smallbreak\hskip 34pt\llap{#1}\kern 5pt}
\def\puce{\smallbreak\hskip 6pt{$\scriptstyle\bullet$}\kern 5pt}

\def\eqalign#1{\null\,\vcenter{\openup\jot\m@th\ialign{
\strut\hfil$\displaystyle{##}$&$\displaystyle{{}##}$\hfil
&&\quad\strut\hfil$\displaystyle{##}$&$\displaystyle{{}##}$\hfil
\crcr#1\crcr}}\,}

\catcode`\@=12

\showboxbreadth=-1  \showboxdepth=-1

\newcount\numerodesection \numerodesection=1
\def\section#1{\bigbreak
 {\bf\number\numerodesection.\ \ #1}\nobreak\medskip
 \advance\numerodesection by1}

\mathcode`A="7041 \mathcode`B="7042 \mathcode`C="7043 \mathcode`D="7044
\mathcode`E="7045 \mathcode`F="7046 \mathcode`G="7047 \mathcode`H="7048
\mathcode`I="7049 \mathcode`J="704A \mathcode`K="704B \mathcode`L="704C
\mathcode`M="704D \mathcode`N="704E \mathcode`O="704F \mathcode`P="7050
\mathcode`Q="7051 \mathcode`R="7052 \mathcode`S="7053 \mathcode`T="7054
\mathcode`U="7055 \mathcode`V="7056 \mathcode`W="7057 \mathcode`X="7058
\mathcode`Y="7059 \mathcode`Z="705A

% handling accented characters in plain TeX :

\def\diagram#1{\def\normalbaselines{\baselineskip=0pt\lineskip=5pt}
\matrix{#1}}

\def\vfl#1#2#3{\llap{$\textstyle #1$}
\left\downarrow\vbox to#3{}\right.\rlap{$\textstyle #2$}}

\def\hfl#1#2#3{\smash{\mathop{\hbox to#3{\rightarrowfill}}\limits
^{\textstyle#1}_{\textstyle#2}}}

\def\ogoth{{\goth o}}

\def\pgoth{{\goth p}}

\def\Q{{\bf Q}}
\def\Qp{\Q_p}

\def\N{{\bf N}}

\def\Z{{\bf Z}}
\def\Zp{\Z_p}
\def\F{{\bf F}}
\def\Fp{{\F_{\!p}}}

\def\Hom{\mathop{\rm Hom}\nolimits}

\def\Card{\mathop{\rm Card}\nolimits}
\def\Gal{\mathop{\rm Gal}\nolimits}

\def\droite#1{\,\hfl{#1}{}{8mm}\,}

\def\to{\rightarrow}

\def\normressym(#1,#2)_#3{\displaystyle\left({#1,#2\over#3}\right)}

\def\mod{\mathop{\rm mod.}\nolimits}
\def\pmod#1{\;(\mod#1)}

\newcount\refno 
\long\def\ref#1:#2<#3>{                                        
\global\advance\refno by1\par\noindent                              
\llap{[{\bf\number\refno}]\ }{#1} \pointir{\it #2} #3\goodbreak }

\def\citer#1(#2){[{\bf\number#1}\if#2\empty\relax\else,\ {#2}\fi]}

\newbox\bibbox
\setbox\bibbox\vbox{\bigbreak
\centerline{{\pc BIBLIOGRAPHIC} {\pc REFERENCES}}

\ref{\pc CASSELS} (J):
Local fields,
<Cambridge University Press, Cambridge, 1986. xiv+360 pp.>
\newcount\cassels \global\cassels=\refno

\ref{\pc DALAWAT} (C):
Local discriminants, kummerian extensions, and elliptic curves,
<Journal of the Ramanujan Mathematical Society, {\bf 25} (2010) 1,
pp.~25--80. Cf.~arXiv\string:0711.3878.>      
\newcount\locdisc \global\locdisc=\refno

\ref{\pc DALAWAT} (C):
Further remarks on local discriminants, 
<0909.2541.>    
\newcount\further \global\further=\refno

\ref{\pc DALAWAT} (C):
Serre's ``{\it formule de masse\/''} in prime degree,
<1005.2016.>    
\newcount\csdmass \global\csdmass=\refno

\ref{\pc DEL \pc CORSO} (I) and {\pc DVORNICICH} (R):
The compositum of wild extensions of local fields of prime degree,
<Monatsh.\ Math.\ {\bf 150} (2007) 4, pp.~271--288.>
\newcount\delcorso \global\delcorso=\refno

\ref{\pc NEUMANN} (O):
Two proofs of the Kronecker-Weber theorem ``according to Kronecker, and
Weber'',
<J.\ Reine Angew.\ Math.\  {\bf 323}  (1981), 105--126.>
\newcount\neumann \global\neumann=\refno

\ref{\pc SERRE} (J-P):
Corps locaux,
<Publications de l'Universit{\'e} de Nancago, No.~{\sevenrm VIII}, Hermann,
Paris, 1968, 245 pp.>
\newcount\corpslocaux \global\corpslocaux=\refno

\ref{\pc SERRE} (J-P):
Une ``formule de masse" pour les extensions totalement ramifi{\'e}es de
degr{\'e} donn{\'e} d'un corps local, 
<Comptes Rendus {\bf 286} (1978), pp.~1031--1036.>
\newcount\serremass \global\serremass=\refno

} %\bibbox

\centerline{\bf Final remarks on local discriminants} 
\bigskip\bigskip 
\centerline{Chandan Singh Dalawat} 
\centerline{\it Harish-Chandra Research Institute}
\centerline{\it Chhatnag Road, Jhunsi, Allahabad 211019, India} 
\centerline{\it dalawat@gmail.com}

\bigskip\bigskip

{{\bf Abstract}. We show how the ramification filtration on the maximal
  elementary abelian $p$-extension ($p$ prime) on a local number field of
  residual characteristic~$p$ can be derived using only Kummer theory and a
  certain orthogonality relation for the Kummer pairing, even in the absence
  of a primitive $p$-th root of~$1$~; the case of other local fields was
  treated earlier.  In all cases, we compute the contribution of cyclic
  extensions to Serre's degree-$p$ mass formula.  \footnote{}{Keywords~: Local
    fields, elementary abelian $p$-extensions, ramification filtration,
    discriminants, Serre's mass formula.}}

\bigskip

{\bf 1. Introduction}\pointir Let $p$ be a prime number, and $K$ a local
number field or a local function field of residual characteristic~$p$, so that
$K$ is a finite extension of $\Qp$ or of $\Fp(\!(\pi)\!)$, where $\pi$ is
transcendental.  Let $M$ be the maximal elementary abelian $p$-extension of
$K$, and $G=\Gal(M|K)$.  The profinite group $G$ comes with a natural
filtration $(G^u)_{u\in[-1,+\infty[}$ (in the upper numbering).  Local class
field theory provides an isomorphism $K^\times\!/K^{\times p}\to G$ preserving
the filtrations and thereby determines $(G^u)_u$.

But a more elementary derivation is possible.  Namely, when $K$ is a local
function field, a certain orthogonality relation for the Artin-Schreier
pairing allows us to determine the filtration $(G^u)_{u}$ in terms of the
filtration on $K/\wp(K)$ \citer\further(), \S5.  Also, when $K$ is a local
number field containing a primitive $p$-th root $\zeta$ of~$1$, the analogous
orthogonality relation for the Kummer pairing allows us to determine
$(G^u)_{u}$ in terms of the filtration on $K^\times\!/K^{\times p}$
\citer\further(), \S4.  

The first purpose of this Note is to determine $(G^u)_{u}$ when $K$ is a local
number field but $\zeta\notin K$.  The idea is to determine the filtered
subspace $V\subset K(\zeta)^\times\!/K(\zeta)^{\times p}$ corresponding to the
exponent-$p$ kummerian extension $M(\zeta)|K(\zeta)$, and then use the
orthogonality relation for the Kummer pairing $G\times V\to{}_p\mu$ to deduce
$(G^u)_{u}$.

The subspace $V$ is determined in~\S2 in a purely algebraic context.  The
filtration on $V$, the orthogonality relation, and the filtration $(G^u)_{u}$
are derived in~\S4.  An overall summary is provided in \S5 to bring out the
analogy between the three cases.

We then compute (\S6) the contribution of degree-$p$ cyclic extensions to
Serre's mass formula for separable degree-$p$ extensions of $K$.  It turns out
that an extension of the ideas in \S2 and \S4 from $K(\zeta)|K$ to
$K(\!\root{p-1}\of{K^\times})|K$, inspired by the paper \citer\delcorso(),
also leads to an elementary proof of the said mass formula, based only on
Kummer theory \citer\locdisc() or Artin-Schreier theory \citer\further() and
some purely algebraic ingredients.  See \citer\csdmass().

\bigskip

{\bf 2.  Algebraic preliminaries}\pointir We shall need some purely algberaic
results which are often used in the proof of the local \citer\cassels(p.~155)
or the global \citer\neumann(p.~110) Kronecker-Weber theorem.  Our
presentation is intrinsic, and shows the equivalence of the statements in
these two sources.

Let $p$ be a prime number and let $F$ be a field in which $p$ is invertible.
We want to understand the degree-$p$ cyclic extensions of $F$ in terms of
those of $K=F(\zeta)$, where $\zeta$ is a primitive $p$-th root of~$1$.

Let $E|F$ be a cyclic extension of degree $p$.  The extension
$E(\zeta)|K$ is also cyclic of degree~$p$~; it corresponds therefore
(``Kummer theory'') to an $\Fp$-line $D\subset
K^\times\!/K^{\times p}$.  The group $\Delta=\Gal(K|F)$
acts on the latter space.  Let $\omega:\Delta\to\F_p^\times$ be the cyclotomic
character giving the action of $\Delta$ on the $p$-th roots of $1$, so that
$\sigma(\zeta)=\zeta^{\omega(\sigma)}$ for every $\sigma\in\Delta$.  

Sometimes we think of the target of $\omega$ as being the interval
$[1,p[\;\,\subset\Z$.  For $p$ odd, $(\Z/p\Z)^\times$ is often identified with
the torsion subgroup of $\Zp^\times$.

\th LEMMA 1
\enonce
The\/ $\Fp$-line\/ $D$ is\/ $\Delta$-stable, and\/ $\Delta$ acts on\/ $D$
via\/ $\omega$. 
\endth
Let $ a\in K^\times$ be such that its image $\bar a$ modulo
$K^{\times p}$ generates $D$, and $x$ a $p$-root of $a$, so that
$E(\zeta)=K(x)$.  We have $\sigma(\bar a)=\overline{\sigma( a)}$ for
every $\sigma\in\Delta$.  We have to first show that $\sigma(\bar a)\in D$.
Identify $\Delta$ with $\Gal(E(\zeta)|E)$.  For every $\sigma\in\Delta$, we
have  
$$
K(x)=
F(\zeta,x)=
F(\sigma(\zeta),\sigma(x))=
K(\sigma(x))
$$
and $(\sigma(x))^p=\sigma(x^p)=\sigma( a)$ is in $K^\times$,
so $\bar a$ and $\overline{\sigma( a)}$ belong to the same
$\Fp$-line, namely $D$.  Hence $D$ is\/
$\Delta$-stable. 

Let $\eta:\Delta\to\F_p^\times$ be the character through which $\Delta$ acts
on $D$, so that, for a generator $\tau$ of $\Delta$, we have
$\tau(x)=bx^{\eta(\tau)}$ for some $b\in K^\times$.  Let $g$ be the generator
$x\mapsto\zeta x$ of the group $\Gal(K(x)|K)$, so that $g(\zeta)=\zeta$.
Hence $\tau(g(x))=\tau(\zeta x)=\zeta^{\omega(\tau)}bx^{\eta(\tau)}$, on the
one hand.

On the other hand, $g(\tau(x)) =g(bx^{\eta(\tau)})
=b\zeta^{\eta(\tau)}x^{\eta(\tau)}$.  But $\tau g=g\tau$~; comparing the two
computations, we get $\eta=\omega$.  Cf.~\citer\cassels(p.~155).

Conversely,
\th LEMMA 2
\enonce
For every\/ $\Delta$-stable\/ $\Fp$-line $D\subset
K^\times\!/K^{\times p}$ on which $\Delta$ acts via $\omega$,
there is a (unique) degree-$p$ cyclic extension $E|F$ such that
$K(\!\root p\of D)=E(\zeta)$. 
\endth
Keep the notation of lemma~1.  Notice first that the extension $K(x)$ is
galoisian over $F$, for it contains, for every $\sigma\in\Delta$, a $p$-th
root of $\sigma(a)$, namely $bx^{\omega(\sigma)}$, where $b\in K^\times$ is
such that $\sigma(a)=a^{\omega(\sigma)}b^p$.  If $\Gal(K(x)|F)$ is commutative
(it would then have to be cyclic because the orders of $G$ and $\Delta$ are
relatively prime), the fixed field under the index-$p$ subgroup would be the
desired $E$.

Let us show that $\Gal(K(x)|F)$ is indeed commutative.  We have seen that
$\tau(x)=bx^{\omega(\tau)}$ (for some $b\in K^\times$) and $g(x)=\zeta x$ (for
some $\zeta\in{}_p\mu$).  Therefore $\tau g(x)=\tau(\zeta
x)=\zeta^{\omega(\tau)}bx^{\omega(\tau)}$.  On the other hand,
$g\tau(x)=g(bx^{\omega(\tau)})=b\zeta^{\omega(\tau)}x^{\omega(\tau)}$.  So the
extension $K(x)|F$ is abelian, as claimed.

\smallskip

\th COROLLARY 3
\enonce
Let\/ $V$ be the\/ $\omega$-eigenspace for the action of\/ $\Delta$ on\/
$K^\times\!/K^{\times p}$.  The map\/ $E\mapsto D$ defined by\/
$E(\zeta)=K(\!\root p\of D)$ is a bijection of the set of degree-$p$
cyclic extensions\/ $E|F$ onto the set of\/ $\Fp$-lines\/ $D\subset V$.
\endth

As usual, $V$ can also be written as the image of a certain projector
$\varepsilon\in\Fp[\Delta]$, when $K^\times\!/K^{\times p}$ is regarded as an
$\Fp[\Delta]$-module.  Indeed, let $\varepsilon
=(1/m)\sum_{\sigma\in\Delta}\omega(\sigma^{-1})\sigma\in\Fp[\Delta]$, where
$m$ is the order of $\Delta$, so that $m$ divides $p-1$.  It is easily
verified that $\tau\varepsilon=\omega(\tau)\varepsilon$ for every
$\tau\in\Delta$~:
$$
\eqalign{
\tau\varepsilon
&={1\over m}\sum_{\sigma\in\Delta}\tau.\omega(\sigma^{-1})\sigma\cr 
&={1\over m}\sum_{\sigma\in\Delta}\omega(\sigma^{-1})\tau\sigma\cr 
&={\omega(\tau)\over m}
   \sum_{\sigma\in\Delta}\omega((\tau\sigma)^{-1})\tau\sigma\cr  
&={\omega(\tau)\over m}m\varepsilon=\omega(\tau)\varepsilon.\cr
}$$
Upon multiplying both sides by $\omega(\tau^{-1})$ and summing over
$\tau\in\Delta$, we get
$\varepsilon\sum_{\tau\in\Delta}\omega(\tau^{-1})\tau=m\varepsilon$, or
$\varepsilon.m\varepsilon=m\varepsilon$, and, as $m$ is invertible in $\Fp$,
we get $\varepsilon^2=\varepsilon$, showing that $\varepsilon$ is an
idempotent.  (It is clear that $\varepsilon(\bar x)=\bar 1$ for every $x\in
F^\times$.) 

In view of this, the preceding lemmas can be reformulated as follows, to bring
out their equivalence with \citer\neumann(p.~110).

\th LEMMA 4
\enonce
For every degree-$p$ cyclic extension\/ $E|F$, there is a unique\/ $\Fp$-line 
$D\subset\varepsilon(K^\times\!/K^{\times p})$ such that\/
$E(\zeta)=K(\!\root p\of D)$.  Conversely, for every\/ $\Fp$-line\/
$D\subset\varepsilon(K^\times\!/K^{\times p})$, there is a
unique degree-$p$ cyclic extension\/ $E|F$ such that\/
$K(\!\root p\of D)=E(\zeta)$.
\endth
[The map $E\mapsto D$ defined by $E(\zeta)=K(\!\root p\of
D)$ gives a bijection between the degree-$p$ cyclic extensions of $F$ and
$\Fp$-lines in $\varepsilon(K^\times\!/K^{\times p})$.]

\smallbreak

Let us summarise.  Let $N$ be the maximal elementary abelian $p$-extension of
$F$.  The results of this~\S\ say that $NK=K(\root p\of V)$, where
$V=\varepsilon(K^\times\!/K^{\times p})$ is the $\omega$-eigenspace for the
action of $\Delta$ on $K^\times\!/K^{\times p}$.  

\bigbreak

{\bf 3. Notations}\pointir The notation to be used in \S4 has been collected
here for reference, and will also be recalled as needed.

--- $F,\zeta,K,\Delta,\omega,\varepsilon,V$~: $F$ is a finite extension of
$\Qp$ which is regular in the sense of Shafarevich~: $\zeta\notin F$~; so
$p\neq2$.  Put $K=F(\zeta)$, and keep the notation $\Delta=\Gal(K|F)$,
$\omega:\Delta\to(\Z/p\Z)^\times$, $\varepsilon\in\Fp[\Delta]$,
$V=\varepsilon(K^\times\!/K^{\times p})$ from \S2.

--- $e,f, e_1$~: $e$ is the ramification index of $F|\Qp$, and $f$ the
residual degree.  We put $e_1=e/(p-1)$, which need not be an integer.

--- $m,s,r$~: $s$ is the ramification index of $K|F$, and $r$ the residual
degree, so $m=sr$ equals $\Card\Delta=[K:F]$.  In particular, $r$ and $s$ are
prime to~$p$.  Notice that the ramification index of $K|\Qp$ is $es$ and the
residual degree $fr$.  Also, $e_1s$ is an integer because $K$ contains
$\zeta$.

--- $k,q,F_r,k_r$~: $k$ is the residue field of $F$, $q=p^f$ is its
cardinality, $F_r$ is the maximal unramified extension of $F$ in $K$, and
$k_r$ is the common residue field of $F_r$ and $K$.  We have $[F_r:F]=r$,
$[K:F_r]=s$.

Let us indicate how $s$ and $r$ can be computed.  Adjoining $\zeta=\root
p\of1$ to $F$ is the same as adjoining $\root{p-1}\of{-p}$
\citer\locdisc(prop.~24), so that the degree $[K:F]=m$ equals the order of
$\overline{-p}\in F^\times\!/F^{\times p-1}$~; this order divides $p-1$.
Momentarily let $s>0$ be the smallest integer such that $\overline{-p}^s\in
k^\times\!/k^{\times p-1}$~; it is also the smallest integer such that
$p-1\,|\,es$.  It is clear that the ramification index of $K|F$ is $s$ and the
residual degree $r=m/s$.  Indeed, if $u\in F^\times$ is a unit such that $\bar
u=\overline{-p}^s$, then the order of $\bar u\in k^\times\!/k^{\times p-1}$ is
$r$, and $K$ contains $\root{p-1}\of u$, in the shape of
$(\!\root{p-1}\of{-p})^s$.  On the one hand, the extension
$F_r=F(\!\root{p-1}\of u)$ is unramified of degree~$r$ over $F$.  On the other
hand, the extension $K|F_r$ is totally ramified of degree~$s$ because
$\overline{-p}\in F_r^\times\!/F_r^{\times p-1}$ has order $s$, and $s$ is the
smallest exponent ($>0$) such that $\overline{-p}^s\in
k_r^\times\!/k_r^{\times p-1}$.

--- $\varpi,\pi,v, w$~: $\varpi$ is a uniformiser of $F$, $\pi$ is the
uniformiser $1-\zeta$ of $\Qp(\zeta)$, $v$ is the valuation on $F$ such that
$v(\varpi)=1$ and $w$ is the valuation on $K$ such that $w(\varpi)=s$.

--- $\ogoth,{\goth O},\pgoth,{\goth P},U_i,\bar U_i$~: Let $\ogoth,{\goth O}$
be the rings of integers of $F,K$ and $\pgoth,{\goth P}$ their unique maximal
ideals.  The $\Fp$-space $\bar U_0=K^\times\!/K^{\times p}$ comes with the
filtration $(\bar U_i)_{i>0}$, where $U_i=1+{\goth P}^i$~; we also have
$\overline{{\goth O}^\times}=\bar U_1$.  We have $\bar U_i=\{1\}$ for
$i>pe_1s$, and $\bar U_{pe_1s}$ is a line \citer\locdisc(prop.~42).

--- $\mu,n,\xi,\bar\mu$~: $\mu$ is the $p^\infty$-torsion subgroup of ${\goth
  O}^\times$~; it is cyclic of order $p^n$, and $\xi$ is a generator. We have
$n>0$ by hypothesis, but $n>1$ is possible.  (As an example, consider
$F=\Qp(\xi)^\Delta$, where $\Delta=(\Z/p\Z)^\times$ is identified as a
subgroup of $\Gal(\Qp(\xi)|\Qp)=(\Z/p^m\Z)^\times$ and $\xi$ is a primitive
$p^m$-th root of~$1$.  As $[F:\Qp]=p^{m-1}$, it does not contain $\zeta$, but
$F(\zeta)=\Qp(\xi)$, so $n=m$).  Finally, $\bar\mu\subset\bar U_1$ is the
image of $\mu$~; it is generated by $\bar\xi$.

--- $a(i),b^{(i)}$~: For every $i>0$, put $\displaystyle
a(i)=\left\lfloor{i-1\over p-1}\right\rfloor$ and $b^{(i)}=i+a(i)$.  Notice
that $i\mapsto b^{(i)}$ is an increasing bijection of $\N^*$ with the set of
integers in $\N^*$ prime to~$p$.

--- $N,G,H$~: $N$ is the maximal elementary abelian $p$-extension of $F$,
    $G=\Gal(N|F)$, $H=\Gal(NK|K)$.

\bigbreak

{\bf 4. Regular local number fields}\pointir Recall that lines in $V\subset
K^\times\!/K^{\times p}$ correspond to degree-$p$ cyclic extensions of $F$.
We would like to determine the filtration $V_i=V\cap\bar U_i$ on $V$.  

We shall go about it slowly, treating some special cases first in order to
bring out the essential ideas.  Almost the only thing we use is the fact that
the unique ramification break of a degree-$p$ cyclic extension of $F$ occurs
at $-1$ or at $b^{(i)}$ for some $i\in[1,e]$.

With the notation introduced in \S2, it is clear that $V_{pi}=V_{pi+1}$ for
every $i\neq e_1s$ and that $V_{pe_1s+1}=\{\bar1\}$ \citer\locdisc(prop.~42).
Let us show first that $V\subset\bar U_1$ and that it contains the lines
$\bar\mu,\bar U_{pe_1s}$.

\th LEMMA 5
\enonce
With these notations, $\bar U_{pe_1s}\subset V\subset\bar U_1$ and\/ $\bar
\mu\subset V$.
\endth

$\bar\mu\subset V$~: Let $\tau\in\Delta$ be a generator and let
$g\in(\Z/p^n\Z)^\times$ be such that $\tau(\xi)=\xi^g$~; we have to show that
$g\equiv\omega(\tau)\pmod p$.  Now, $\zeta=\xi^{p^{n-1}}$ is a primitive $p$-th
root of~$1$, so $\tau(\zeta)=\zeta^{\omega(\tau)}=\xi^{\omega(\tau)p^{n-1}}$, by
hypothesis.  But we also have
$\tau(\zeta)=\tau(\xi)^{p^{n-1}}=\xi^{gp^{n-1}}$, so that
$gp^{n-1}\equiv\omega(\tau)p^{n-1}\pmod{p^n}$, and hence $g\equiv\omega(\tau)\pmod
p$.

$V\subset\bar U_1$~: Let $w:K^\times\to\Z$ be the surjective valuation, so
that $\bar x\in\bar U_1\Leftrightarrow w(x)\in p\Z$ ($x\in K^\times$).
Suppose that $\bar x\in V$ for some $x\in K^\times$, and let $\tau\in\Delta$
be a generator.  By hypothesis, $\tau(x)=x^{\omega(\tau)}y^p$ for some $y\in
K^\times$.  Taking valuations, we get $w(x)\equiv w(x)\omega(\tau)\pmod p$,
and, as $\omega(\tau)\not\equiv1\pmod p$, we conclude that $w(x)\equiv0\pmod
p$, and hence $\bar x\in\bar U_1$.

$\bar U_{pe_1s}\subset V$~: Let $\alpha\in\ogoth^\times$ be a unit of $F$ the
trace of whose reduction $S_{k|\Fp}(\hat\alpha)\neq0$ in $\Fp$.  When $\alpha$
is considered as a unit of $K$, we have
$S_{k_r|\Fp}(\hat\alpha)=rS_{k|\Fp}(\hat\alpha)$, which is still $\neq0$, for
$r\not\equiv0\pmod p$.  The line $\bar U_{pe_1s}$ is therefore generated by
the image of $1+\alpha p\pi$, where $\pi=1-\zeta$~; see the discussion after
\citer\locdisc(prop.~33).  But for every $\tau\in\Delta$, we have
$\tau(\pi)\equiv\omega(\tau)\pi\pmod{\pi^2}$, so
$$
\tau(1+\alpha p\pi)
\equiv1+\omega(\tau)\alpha p\pi 
\equiv(1+\alpha p\pi)^{\omega(\tau)}
\pmod{{\goth P}^{pe_1s+1}}
$$
(recalling that $\pi{\goth O}={\goth P}^{e_1s}$, $p{\goth O}={\goth P}^{es}$,
$e=(p-1)e_1$), which shows that $\tau(1+\alpha p\pi)=(1+\alpha
p\pi)^{\omega(\tau)}\beta^p$ for some $\beta\in U_1$, and hence
$\overline{1+\alpha p\pi}\in V$.

\smallskip

Consider for a moment the special case $F=\Qp$, so that $e_1=1$ and
$\pi=1-\zeta$ is a uniformiser of $K$.  It is easily seen that
$V_p=V_{p-1}=\cdots=V_2$~: if $\bar x\in V_i$ for some $i\in[2,p]$ and some
$x\equiv1+\alpha\pi^j\pmod{\pi^{j+1}}$ ($\alpha\in\Zp^\times$), then computing
$\tau(\bar x)$ in two different ways leads to the result.

Indeed, working $\pmod{\pi^{j+1}}$, we have
$\tau(x)\equiv1+\alpha\tau(\pi)\equiv1+\alpha\omega(\tau)^j\pi^j$ on the one
hand.  On the other hand, as $\bar x\in V$, we have
$\tau(x)=x^{\omega(\tau)}y^p$ for some $y\in U_1$.  But $y^p\equiv1$, so
$\tau(x)\equiv(1+\alpha\pi)^{\omega(\tau)}\equiv1+\omega(\tau)\alpha\pi$.  The two
computations imply $\omega(\tau)^{j-1}\equiv1\pmod\pi$, and in fact $\!\!\pmod
p$, because $\omega(\tau)\in\Zp^\times$.  But $\omega(\tau)^{j-1}\equiv1\pmod p$
holds for a generator $\tau\in\Delta$ (which has order $p-1$) only when $j=p$.
It follows that $V_p=V_2$ and hence $V=\bar U_p\bar\mu$
\citer\cassels(p.~156)~; the line $\bar U_p$ corresponds to the unramified
degree-$p$ extension, and the line $\bar\mu$ to the cyclotomic
$(\Z/p\Z)$-extension.

The result $V_2=\bar U_p$ (when $F=\Qp$) could also have been obtained by
remarking that in this case the unique ramification break of a ramified
degree-$p$ cyclic extension $L|K$ coming from $F$ occurs at~$p-1$
\citer\locdisc(prop.~63), an argument which works for any finite extension
$F|\Qp$ (such that $\zeta\notin F$).  Hence the following bit of information
about the filtration on~$V$~:

\th LEMMA 6
\enonce
We have\/ $\bar U_{pe_1s}=V_{pe_1s-1}=\cdots=V_{pe_1s-s+1}$, where\/ $s$ is
the ramification index of\/ $K|F$.    
\endth
Let $D\subset V_{pe_1s-s+1}$ be a line such that $D\neq\bar U_{pe_1s}$, $E|F$
the corresponding degree-$p$ cyclic extension, and $t>0$ its unique
ramification break.  The unique ramification break of $E(\zeta)|K$
occurs at $ts$ \citer\locdisc(proof of prop.~63) on the one hand, and at $s-i$
for some  $i\in[1,s[$, on the other \citer\locdisc(prop.~60).  But $ts=s-i$ is
impossible, so there is no such $D$, and hence $\bar U_{pe_1s}=V_{pe_1s-s+1}$.

Let us next determine the $\Fp$-dimension of $V_{pe_1s-s}/V_{pe_1s-s+1}$.

\th LEMMA 7
\enonce
We have\/ $\dim_{\Fp}V_{pe_1s-s}/V_{pe_1s-s+1}=f$, where\/ $f=[k:\Fp]$ is the
residual degree of\/ $F|\Qp$.
\endth
Let $\varpi$ be a uniformiser of $F$ and recall that $\pi=1-\zeta$, where
$\zeta\in K$ is a primitive $p$-th root of~$1$.  We have
$w(p\pi\varpi^{-1})=pe_1s-s$, so that  $p\pi\varpi^{-1}$ is an
$\goth O$-basis of ${\goth P}^{pe_1s-s}$.
 
Use this basis to identify $\bar U_{pe_1s-s}/\bar U_{pe_1s-s+1}$ with
$k_r={\goth O}/{\goth P}$, the residue field of $K$, by sending
$\overline{1+\alpha p\pi\varpi^{-1}}$ ($\alpha\in{\goth O}$) to $\hat\alpha$.
We claim that then $V_{pe_1s-s}/V_{pe_1s-s+1}$ gets identified with the
subspace $k\subset k_r$.  The idea is to show that the $\Fp[\Delta]$-module
structure on $k_r$ coming from the said identification is the usual structure
twisted by $\omega$.

Let $\tau\in\Delta$ be a generator.  Recall that
$\tau(\pi)\equiv\omega(\tau)\pi\pmod{\pi^2}$ and that $\tau(\alpha)\equiv
\alpha^q\pmod{\varpi}$ for every integer $\alpha\in F_r$, where $q=\Card k$ is
the residual cardinality of $F$.  Because $pe_1s-s$ is prime to~$p$, the group
$\bar U_{pe_1s-s}/\bar U_{pe_1s-s+1}$ is canonically isomorphic to $
U_{pe_1s-s}/U_{pe_1s-s+1}$ \citer\locdisc(proof of prop.~42) and because $K$
and $F_r$ have the same residue field, every element of
$U_{pe_1s-s}/U_{pe_1s-s+1}$ is represented by $1+\alpha p\pi\varpi^{-1}$ for
some integer $\alpha\in F_r$.  

Now, $\tau(\alpha p\pi\varpi^{-1}) \equiv
\omega(\tau)\alpha^qp\pi\varpi^{-1}\pmod{{\goth P}^{pe_1s-s+1}}$, from which it
follows that $\tau(\hat\alpha)=\omega(\tau){\hat\alpha}^q$ for every
$\hat\alpha\in k_r$.  The $\omega$-eigenspace for this new $\Delta$-action on
$k_r$ is thus $k$, the set of $\hat\alpha\in k_r$ such that
${\hat\alpha}^q=\hat\alpha$.

But $V_{pe_1s-s}/V_{pe_1s-s+1}\subset \bar U_{pe_1s-s}/\bar U_{pe_1s-s+1}$ is
the $\omega$-eigenspace, hence it gets identified with $k\subset k_r$, proving
the lemma.

\medskip

We shall need the following  elementary fact.  

\th LEMMA 8
\enonce
The number of prime-to-$p$ integers in\/ $[1,pe_1[$ is\/ $e$, and they are 
$$
1=b^{(1)}<b^{(2)}<\cdots<b^{(e)}\qquad(<pe_1),
$$
where\/ $b^{(i)}=i+a(i)$ and\/ $\displaystyle a(i)=\left\lfloor{i-1\over
    p-1}\right\rfloor$ for every integer\/ $i\in[1,e]$.
\endth
Write $e=(p-1)c+c'$, with $c\in\N$ and the remainder $c'\in[0,p-1[$, so that 
$e_1=c+c'/(p-1)$, where $c'/(p-1)\in[0,1[$ is rational, and 
$$
pe_1=pc+{pc'\over p-1}=pc+c'+{c'\over p-1}.
$$
It is now clear that the number of integers in\/ $[1,pe_1[$ which are prime
to~$p$ is\/ $pc+c'-c=e$.  That they are precisely $b^{(1)},\ldots,b^{(e)}$ is
left as an exercise.

Recall that the unique ramification break of a ramified degree-$p$ cyclic
extension $L|F$ occurs at $b^{(i)}$ for some $i\in[1,e]$~; see for example
\citer\locdisc(prop.~63).

\th LEMMA 9
\enonce
For\/ $i\in[1,e]$, we have
$\dim_{\Fp}V_{pe_1s-sb^{(i)}}/V_{pe_1s-sb^{(i)}+1}=f$, and\/
$V_{pe_1s-b^{(i-1)}s} = V_{pe_1s-b^{(i)}s+1}$  (with the convention
$b^{(0)}=0$).     
\endth
We have already seen the case $i=1$ of the first part in lemma~7, whose proof
can be adapted to the case $i>1$ by using the ${\goth O}$-basis
$p\pi\varpi^{-b^{(i)}}$ of ${\goth P}^{pe_1s-b^{(i)}s}$ to identify $\bar
U_{pe_1s-b^{(i)}s}/\bar U_{pe_1s-b^{(i)}s+1}$ with $k_r(1)$, the
$\Delta$-module $k_r$ with the action twisted by $\omega$.

The case $i=1$ of the second part is lemma~6, and the same proof works for
$i>1$.  Indeed, let $D\subset V_{pe_1s-b^{(i)}s+1}$ be an $\Fp$-line and
suppose that the unique ramification break of the corresponding degree-$p$
cyclic extension $E|F$ occurs at some $t<b^{(i)}$.  As $t$ is prime to~$p$, we
have $t\le b^{(i-1)}$, and it follows that $D\subset V_{pe_1s-b^{(i-1)}s}$.

\th LEMMA 10
\enonce
We have\/ $V_{pe_1s-b^{(e)}s}=V$.  Equivalently, $V\subset\bar
U_{pe_1s-b^{(e)}s}$. 
\endth
The idea of the proof is the same as for the last few lemmas.  Explicitly, let
$D\subset V$ be an $\Fp$-line, $E|F$ the corresponding degree-$p$ cyclic
extension, and $t$ the unique ramification break of $\Gal(E|F)$~; we have
$t\le b^{(e)}$.  The unique ramification break of $\Gal(E(\zeta)|K)$
occurs at $ts\le b^{(e)}s$, hence $D\subset V_{pe_1s-b^{(e)}s}$.

(We know that $\bar\mu\subset V$ (lemma~5), so lemma~9 leads to the unexpected
consequence that $\bar\mu\subset\bar U_{pe_1s-b^{(e)}s}$.)

\medbreak

Let us pause for a moment to summarise what we have learnt about the
filtration on $V$.  Let us agree to write $B\subset_cA$ if $A$ is an
$\Fp$-space and $B$ is a subspace of $A$ of codimension $c$.  The filtration
on $V$ begins with
$$
\{1\}\subset_1V_{pe_1s}=V_{pe_1s-s+1}\subset_f V_{pe_1s-s}
$$ 
and continues, for every integer $i\in[1,e[$, with
$$
V_{pe_1s-b^{(i)}s}=V_{pe_1s-b^{(i+1)}s+1}\subset_f V_{pe_1s-b^{(i+1)}s},
$$ 
to end finally with $V_{pe_1s-b^{(e)}s}=V$.  It follows that the
$\Fp$-dimension of $V$ is $1+ef=1+[F:\Qp]$.  In short, the breaks in the
filtration $(V_j)_{j>0}$ occur at $j=pe_1s$, where the order of the group
drops by a factor of $p$, and at $j=pe_1s-b^{(i)}s$ for every integer
$i\in[1,e]$, where the order drops by a factor of $q=p^f$.

\medbreak

Now let $N$ be the maximal elementary abelian $p$-extension of $F$, so that
$NK=K(\root p\of V)$ (lemma~4).  Let us briefly indicate how to compute the
ramification filtration on $G=\Gal(N|F)$ using the preceding results, without
any appeal to local class field theory.

In view of the the two cases treated earlier (finite extensions of $\Qp$
having a primitive $p$-th root of~$1$ \citer\locdisc(), local function fields
\citer\further()), it is natural to look for an ``\thinspace orthogonality
relation\thinspace'' for the pairing 
$$
G\times V\to{}_p\mu
$$
which sends $(\sigma,\bar x)$ to $\sigma(\root p\of x)/\root p\of x$, after
having made the identification $G\to H$, where $H=\Gal(NK|K)$.  This is the
content of the prop.~11.

For a subspace $E\subset G$, denote by $E^\perp\subset V$ the subspace such
that $N^EK=K(\root p\of{E^\perp})$.  For example, if $T\subset G$ is the
inertia subgroup (so that $N^T$ is the degree-$p$ unramified extension of
$F$), then $T^\perp=V_{pe_1s}$.  If we identify $G$ with $H$ to get a pairing
$G\times V\to{}_p\mu$, then $E^\perp$ is the orthogonal of $E$.  Denote by
$D^\perp\subset G$ the orthogonal of a subspace $D\subset V$.

\th PROPOSITION 11
\enonce
We have\/ $G^u=G^1$ for\/ $u\in\;]-1,1]$, $G^u=\{1\}$ for\/ $u>b^{(e)}$, and 
$$
(G^u)^\perp=V_{pe_1s-\lceil u\rceil s+1}\qquad(u\in[1,b^{(e)}])
$$
under\/ $G\times V\to{}_p\mu$, the  pairing coming from 
the identification\/ $G\to H$.  
\endth
Notice first that $G^u\neq G$ for $u>-1$, because the quotient $\Gal(F_p|F)$
of $G$ has its break at $-1$, where $F_p$ is the unramified degree-$p$
extension of $F$.  It follows that the index of $G^u\subset G^{-1}$ is~$>1$
for $u>-1$.

Let $u\in\;]-1,1]$ and let $E$ be a hyperplane containing $G^u$, so that
$G/E$ is cyclic of order~$p$.  As the filtration on $G/E$ is the quotient of
the filtration on $G$, the ramification break of $G/E$ occurs somewhere $<u$
(because $G^u\subset E$).  But the only degree-$p$ cyclic extension of $F$
whose ramification break is $<1$ is the unramified one.  So $E=
V_{pe_1s}^\perp$ is the only hyperplane containing $G^u$.  This implies that
$G^u=E=G^1=V_{pe_1s}^\perp$.  

Suppose next that $u>b^{(e)}$~; it suffices to show that every hyperplane
$E\subset G$ contains $G^u$.  Now $G^u\subset E$ if and only if the unique
ramification break of $G/E$ occurs somewhere $<u$.  But this is true for
every $E$, because $u>b^{(e)}$.  Hence $G^u=\{1\}$.

It remains to determine $G^{u\perp}$ for $u\in[1,b^{(e)}]$.  Take a line
$D\neq V_{pe_1s}$ in $V$ and denote by $t\neq-1$ be the unique ramification
break of $G/D^\perp$, so that the unique ramification break of $K(\root p\of
D)|K$ occurs at $ts$.  Then
$$
D\subset G^{u\perp}\Leftrightarrow 
(G/D^\perp)^u=0\Leftrightarrow 
t<u \Leftrightarrow 
ts<\lceil u\rceil s\Leftrightarrow 
D\subset V_{pe_1s-\lceil u\rceil s+1}.
$$
As the two subspaces $G^{u\perp}$ and $V_{pe_1s-\lceil u\rceil s+1}$ of $V$
contain the same lines, they are equal.  Note in particular that
$G^{b^{(e)}\perp}=V_{pe_1s-b^{(e-1)}s}$ (lemma~9), which has codimension~$f$
in~$V$ (lemma~10).

\smallbreak

Now it is an easy matter to determine the filtration on $G$, knowing as we do
the filtration on $V$.

\th COROLLARY 12
\enonce
The upper ramification breaks of\/ $(G^u)_{u\in[-1,+\infty[}$ occur at\/ $-1$ 
and at the\/ $b^{(i)}$ for $i\in[1,e]$~; the codimensions are given by
$$
\{1\}\subset_f
G^{b^{(e)}}\subset_f\cdots\subset_f 
G^{b^{(2)}}\subset_f
G^{b^{(1)}}=G^0\subset_1
G^{-1}=G,
$$
where $\subset_f$ means ``codimension $f$''.  In particular, $G^{pj}=G^{pj+1}$
for every\/~$j$. 
\endth
This follows from prop.~11 and our knowledge of the filtration on~$V$
(lemmas~5--10).

\smallbreak

It is also an easy matter to determine the filtration in the lower numbering
on $G$.  We have the following table for the index of $G^u$ in $G^0$ for
$u\in[0,+\infty[$~:
$$
\vbox{\halign{&\hfil$#$\hfil\quad\cr
u\in&[0,b^{(1)}]&]b^{(1)},b^{(2)}]&\cdots
     &]b^{(e-1)},b^{(e)}]&]b^{(e)},+\infty[\cr
\noalign{\vskip-5pt}
\multispan6\hrulefill.\cr
\hbox{\bf(}G^0:G^u\hbox{\bf)}=&1&q&\cdots&q^{e-1}&q^e\cr
}}
$$
The $e$ positive ramification breaks in the lower numbering occur therefore at
$b_{(i)}=\psi_{N|F}(b^{(i)})$ \citer\corpslocaux(p.~74) for $i\in[1,e]$.  As in
\citer\further(), we have
$$
b_{(i)}=(1+q+\cdots+q^{i-1})+(q^{p-1}+\cdots+q^{a(i)(p-1)}).
$$
\vskip-\baselineskip
\th COROLLARY 13
\enonce
The lower ramification breaks of\/ $(G_l)_{l\in[-1,+\infty[}$ occur at\/ $-1$ 
and at the\/ $b_{(i)}$ for $i\in[1,e]$~; the codimensions are given by
$$
\{1\}\subset_f
G_{b_{(e)}}\subset_f\cdots\subset_f 
G_{b_{(2)}}\subset_f
G_{b_{(1)}}\subset_1
G_{-1}=G.
$$
\endth
An application of \citer\further(lemma~2) now gives the exponent $v_N({\goth
  D}_{N|F})$ of the different ${\goth D}_{N|F}$ as 
$$
v_N({\goth D}_{N|F})=(1+b^{(e)})q^e-(1+b_{(e)})
$$
and the exponent $v(d_{N|F})$ of the discriminant as
$v(d_{N|F})=p.v_N({\goth D}_{N|F})$, because the residual degree of $N|F$
is~$p$.  

Local class field theory was needed in \citer\further(after prop.~5) to
compute the filtration on $G$ and thereby obtain these values for the exponent
of the different and the discriminant.

\bigbreak

{\bf 5. Overall summary}\pointir Let $p$ be a prime number, $K$ a finite
extension of $\Qp$ or of $\Fp(\!(\pi)\!)$, $M$ the maximal elementary abelian
$p$-extension of $K$, and $G=\Gal(M|K)$.  We have seen that it is possible to
determine the filtration (in the upper numbering) on $G$ using only Kummer
theory in the local number field case, and only Artin-Schreier theory in the
local function field case.  In the former case --- where $G$ is finite --- the
lower numbering can also be determined.  In the latter case, one can determine
the lower numbering on the finite quotients of~$G$.

\smallbreak

Consider first a finite extension $K|\Qp$ (of ramification index~$e$ and
residual degree~$f$) containing a primitive $p$-root $\zeta$ of~$1$.  The
filtration on $\overline{K^\times}=K^\times\!/K^{\times p}$ is easily
determined and looks like
$$
\{1\}
\subset_1\bar U_{pe_1}
\subset_f\bar U_{b^{(e)}}
\subset_f\cdots
\subset_f\bar U_{b^{(1)}}
\subset_1\bar U_0=\overline{K^\times},
$$
where $e_1=e/(p-1)$ \citer\locdisc(prop.~42).  If a line
$D\subset\overline{K^\times}$ is in $\bar U_m$ but not in $\bar U_{m+1}$
(which forces $m=0$ or $m=b^{(i)}$ for some $i\in[1,e]$ or $m=pe_1$), then the
unique ramification break of $K(\root p\of D)$ occurs at $pe_1-m$ if $m\neq
pe_1$ \citer\locdisc(prop.~60), at $-1$ if $m=pe_1$ \citer\locdisc(prop.~16).
The filtration $(G^u)_u$ is completely determined by $G^u=G^1$ for\/
$u\in\;]-1,1]$, $G^u=\{1\}$ for\/ $u>pe_1$ and the orthogonality relation
$$
(G^u)^\perp=\bar U_{pe_1-\lceil u\rceil+1}
$$
for $u\in[1,pe_1]$ \citer\locdisc(Part~IX), under the Kummer pairing
$G\times\overline{K^\times}\to{}_p\mu$.  This leads to the description
$$
\{1\}\subset_1
G^{pe_1}\subset_f
G^{b^{(e)}}\subset_f\cdots\subset_f 
G^{b^{(2)}}\subset_f
G^{b^{(1)}}\subset_1G.
$$
The ramification breaks in the lower numbering occur at $-1$, at the
$b_{(i)}$ for $i\in[1,e]$, and at $b_{(e)}+q^e$, where $q=p^f$
\citer\further(prop.~3).

\smallbreak

Consider next a finite extension $K|\Fp(\!(\pi)\!)$ (of residual degree~$f$).
The filtration on $\overline{K}=K/\wp(K)$ looks like
$$
\{0\}
\subset_1\overline\ogoth
\subset_f\overline{\pgoth^{-b^{(1)}}}
\subset_f\overline{\pgoth^{-b^{(2)}}}
\subset_f\cdots
\subset\overline{K}
$$
\citer\further(),~\S6.  If a line $D\subset\overline{K}$ is in
$\overline{\pgoth^{-m}}$ but not in $\overline{\pgoth^{-m+1}}$ (which forces
$m=0$ or $m=b^{(i)}$ for some $i\in\N^*$), then the unique ramification break
of $K(\wp^{-1}(D))$ occurs at $m$ if $m\neq0$ \citer\further(prop.~14), at
$-1$ if $m=0$ \citer\further(prop.~12).  The filtration $(G^u)_u$ is
completely determined by $G^u=G^1$ for $u\in\;]-1,1]$ and the orthogonality
relation 
$$
(G^u)^\perp=\overline{\pgoth^{-\lceil u\rceil+1}}
$$
for $u>0$ \citer\further(prop.~{17}), under the Artin-Schreier pairing
$G\times\overline{K}\to\F_p$, leading to the description
$$
\{1\}\subset\cdots\subset_f 
G^{b^{(2)}}\subset_f
G^{b^{(1)}}\subset_1G.
$$
For $m\in\N$, the breaks in the lower numbering on
$K(\wp^{-1}(\pgoth^{-m}))$ occur at $-1$ and at $b_{(i)}$ for $i\in[1,c(m)]$,
where $c(m)=m-\left\lfloor{m/p}\right\rfloor$ \citer\further(prop.~19).

\smallbreak

Consider lastly a finite extension $K|\Qp$ (of ramification index~$e$ and
residual degree~$f$) {\it not\/} containing $\zeta$, such as the $F$ in
\S3--4, and put $L=K(\zeta)$.  We have determined the filtered subspace
$V\subset\overline{L^\times}$ (\S2, \S4) lines $D$ in which correspond to
degree-$p$ cyclic extensions $E$ of $K$ by the rule $L(\root p\of
D)=E(\zeta)$.  Lines in $V$ correspond therefore to hyperplanes in $G$ and
lead to an orthogonality relation (prop.~11) which determines $(G^u)_u$
(cor.~12) and $(G_l)_l$ (cor.~13).

The information carried by $V$ can be succintly expressed by posing
$W_i=V_{pe_1s-b^{(i)}s}$ for $i\in[0,e]$, where $e_1=e/(p-1)$, $s$ is the
ramification index of $L|K$, and $b^{(0)}=0$ by convention.  We then have the
picture
$$
\{1\}\subset_1 W_0\subset_f W_1\subset_f\cdots\subset_f W_e=V;
$$
the line $W_0$ corresponds to the unramified degree-$p$ extension of $K$, and,
for every $i\in[1,e]$ and every line $D\subset W_i$ such that $D\not\subset
W_{i-1}$, the unique ramification break of the corresponding degree-$p$ cyclic
extension $E|K$ occurs at $b^{(i)}$.  In particular,
$v_{K}(d_{E|K})=(p-1)(1+b^{(i)})$. 

\bigbreak

{\bf 6. The contribution of cyclic extensions}\pointir Let $p$ be a prime
number, $k|\Fp$ a finite extension, $q=p^f=\Card k$, and let $F$ be a local
field with residue field $k$.  The preceding considerations can be applied to
computing the contribution of cylic extensions to Serre's degree-$p$ ``mass
formula'' \citer\serremass().

Recall that the formula in question asserts that $\sum_L q^{-c(L)}=n$, where
$L$ runs through totally ramified extensions of $F$ (in a fixed separable
closure) of degree~$n=[L:F]$, and $c(L)=v_F(d_{L|F})-(n-1)$.  One may ask for
the contribution of {\it cyclic\/} extensions to this formula~; the foregoing
summary (\S5) makes it possible to compute it.

[In the case $p=2$, every separable quadratic extension is cyclic, so the
contribution should be 100\%; this has been verified in the characteristic-$0$
case \citer\locdisc(lemma~67). We shall see that in the characteristic-$2$
case it amounts to the identity
$$
{2+2^2+\cdots+2^f\over 2^{(2-1)f}}+
%{2^{f+1}+2^{f+2}+\cdots+2^{2f}\over 2^{f.(2.2-1)}}+
\cdots+
{2^{(i-1)f+1}+2^{(i-1)f+2}+\cdots+2^{if}\over 2^{(2i-1)f}}+
\cdots
=2.]
$$

Consider first the characteristic-$p$ case.

\th PROPOSITION 14
\enonce
Let\/ $F=k(\!(\pi)\!)$.  When\/ $L$ runs through ramified cyclic
extensions of\/~$F$ of degree~$p$, we have
$$
\sum_L q^{-c(L)}
={p\over q}\cdot{q-1\over p-1}\cdot\sum_{i>0}q^{i-(p-1)b^{(i)}},\leqno{(1)}
$$
where\/ $b^{(i)}=i+a(i)$ and\/ $a(i)=\lfloor(i-1)/(p-1)\rfloor$ for every 
integer\/ $i>0$.
\endth
The idea of the proof is clear from \S5.  Each (ramified, degree-$p$, cyclic)
extension $L|F$ has a unique ramification break, which equals $b^{(i)}$ for
some $i>0$~; if
so, then $c(L)=(p-1)b^{(i)}$.  These $L$ correspond to $\Fp$-lines
$D\subset\overline{\pgoth^{-b^{(i)}}}$ such that
$D\not\subset\overline{\pgoth^{-b^{(i-1)}}}$.  As the dimension of
$\overline{\pgoth^{-b^{(i)}}}$ is $1+if$, the number of such lines $D$ is
$pq^{i-1}+p^2q^{i-1}+\cdots+p^fq^{i-1}=pq^{i-1}(q-1)/(p-1)$.  So the
contribution of such $L$ (or such $D$) to the sum is
$$
{pq^{i-1}(q-1)\over(p-1)q^{(p-1)b^{(i)}}}
={p\over q}\cdot{q-1\over p-1}\cdot q^{i-(p-1)b^{(i)}},
$$
and summing over all $i>0$ gives the result.  Note that, when $p=2$,
$i-(p-1)b^{(i)}=i-b^{(i)}=-a(i)=1-i$, so
$\sum_{i>0}q^{i-(p-1)b^{(i)}}=q/(q-1)$. 

\medskip

Consider next the characteristic-$0$ case of a finite extension $F|\Qp$ with
ramification index~$e$ and residual degree~$f$~; put\/ $e_1=e/(p-1)$ and
$q=p^f$.

\th PROPOSITION 15
\enonce
Suppose that\/ $F|\Qp$ contains a primitive\/ $p$-th root of\/~$1$.  When\/
$L$ runs through ramified degree-$p$ cyclic extensions of\/ $F$, we have
$$
\sum_L q^{-c(L)}
={p\over q^{(p-1)e}}
+{p\over q}{q-1\over p-1}\sum_{i\in[1,e]}q^{i-(p-1)b^{(i)}}.\leqno{(2)}
$$
\endth
Ramified cyclic degree-$p$ extensions $L|F$ are of two kinds.  If the unique
ramification break $t$ of $\Gal(L|F)$ is prime to $p$, then $t=b^{(i)}$ for
some $i\in[1,e]$~; they are called {\it peu ramifi{\'e}es}, correspond to
lines in the $\Fp$-space $\bar U_1$ other than the line $\bar U_{pe_1}$, and
contribute the second term on the right in (2), as we saw in the
characteristic-$p$ case (prop.~14).

If, on the other hand, $p|t$, then $t=pe_1$~; such extensions $L|K$ are
called {\it tr{\`e}s ramifi{\'e}es\/} and correspond to $\Fp$-lines
$D\subset K^\times\!/K^{\times p}$ not contained in $\bar U_1$.  The number of
such lines is $pq^e$.  As we then have $c(L)=pe$, this explains the presence
of the first term on the right in (2).

Consider finally the characteristic-$0$ case in the absence of $\root p\of1$.

\th PROPOSITION 16
\enonce
Suppose that\/ $F|\Qp$ does not contain a primitive\/ $p$-th root of\/~$1$.
When\/ $L$ runs through ramified degree-$p$ cyclic extensions of\/ $F$, 
$$
\sum_L q^{-c(L)}
={p\over q}{q-1\over p-1}\sum_{i\in[1,e]}q^{i-(p-1)b^{(i)}}.\leqno{(3)}
$$
\endth 
This follows easily from the last paragraph of \S5 and the proof in
the previous case.  The only difference is that $F$ now has no {\it tr{\`e}s
  ramifi{\'e}es\/} extensions, which explains the absence in (3) of the first
term on the right in (2).  So, in the characteristic-$0$ case, a lesser
proportion of degree-$p$ extensions is cylic if $\root p\of1\notin K$ than if
$\root p\of1\in K$, all other things being equal.

{\it Remark}\pointir Notice that the method allows us to compute the average
$c(L)$ as $L$ runs through ramified degree-$p$ cyclic extensions of $F$ of
some given kind. We illustrate this with the case $F=\Qp(\!\root p\of1)$.  For
every $i\in[1,p]$, there are $p^i$ extensions $L$ such that $c(L)=(p-1)i$, so
the average is
$$
{\sum_{i\in[1,p]}(p-1)ip^i\over\sum_{i\in[1,p]}p^i}=
{p^{p+2}-p^{p+1}-p^p+1\over p^p-1}.
$$
If we want $L$ to be {\it peu ramifi{\'e}\/}, the sums extend only over
$i\in[1,p[$. 

\medbreak
\line{\hfill\hbox{***}\hfill}

\bigbreak
\unvbox\bibbox 

\bye

\th PROPOSITION 18
\enonce 
For every ramified degree-$p$ separable extension\/ $E|F$, there exists a
cyclic extension\/ $F'|F$ of degree dividing\/~$p-1$ such that\/ $EF'|F'$ is
cyclic (of degree\/~$p$) and\/ $EF'|F$ is galoisian.  If\/ $E|F$ is not
cyclic, then it has exactly\/ $p$ conjugates over\/ $F$.
\endth 
The proof of the first part is adapted from \citer\doud().  Let $\tilde E|F$
be the galoisian closure of $E|F$.  Embded $G=\Gal(\tilde E|F)$ in ${\goth
  S}_p$.  If the ramification group $G_1$ has order~$p$, then $G/G_1$ is
cyclic of order dividing $p-1$ by lemma~17 applied to $C=G_1$.  We may then
take $F'=\tilde E^{G_1}$, for we have $\tilde E=EF'$.

Let us therefore show that the $p$-group $G_1$ has order~$p$.  If $G_1$ were
trivial, the extension $\tilde E|F$ would be tame and the ramification index
of $\tilde E|F$ would be prime to~$p$, in contradiction to the fact that the
ramification index of $E|F$ is $p$.  So the order of $G_1$ is divisible by
$p$~; but it is not divisible by $p^2$, for the order $p!$ of ${\goth S}_p$ is
not.  Therefore $G_1$ has order~$p$, and the proof is over.

Finally, we have the exact sequence $1\to G_1\to G\to G/G_1\to1$.  If $E|F$ is
not cyclic, then $[F':F]>1$ and the conjugation action of $G/G_1$ on $G_1$ is
not trivial (if it were trivial, $G$ would be commutative and hence $E|F$
cyclic).  As $G$ is an extension of a cyclic group of order $>1$ dividing
$p-1$ by a cyclic group of order~$p$, it has exactly $p$ subgroups of
index~$p$. ($G$ is the twisted product $G_1\times_\tau(G/G_1)$, where the
character $\tau:G/G_1\to\Fp^\times$ comes from the cojugation action of
$G/G_1$ on $G_1$, so $G$ has at least one subgroup of index of $p$~; it is not
normal, so its orbit under the conjugation action of $G_1$ on index-$p$
subgroups of $G$ consists of $p$ such subgroups.)  It follows that $E|F$ has
exactly $p$ conjugates.

\smallskip

The unique ramification break $t$ ($>0$) of $G_1=\Gal(\tilde E|F')$, and the
orders $p,pe,pef$ of the groups $G_1,G_0,G$ determine the valuation
$v(d_{E|F})$ of the discriminant of $E|F$ by computing $v(d_{\tilde E|F})$
in two different ways, using the {\it Schachtelungssatz\/} along the two
towers $\tilde E|F'|F$, $\tilde E|E|F$.  Indeed, $v_{F'}(d_{\tilde
  E|F'})=(p-1)(1+t)$.  The extension $F'|F$ is tamely ramified of group
$G/G_1$ and inertia subgroup $G_0/G_1$, so the ramification index is $e$ and
the residual degree is $f$, whence $v(d_{F'|F})=(e-1)f$.  The extension
$\tilde E|E$ has the same ramification index and residual degree as $F'|F$, so
$v_E(d_{\tilde E|E})=(e-1)f$.  In conclusion, $v(d_{E|F})=(p-1)(1+t)/e$.
Cf.~\citer\doud(), where $t$ is called $d$ and $e$ is called $t$, and where
$(1+d/t)$ should be interpreted as $(1+d)/t$.

In principle, it should now be possible to prove Serre's degree-$p$ mass
formula by computing the contribution of each such $F'$~; when $F'=F$, then
$E|F$ is cyclic, and the contribution of $F$ itself is given in prop.~14--16.
The number of $F'$ can de deduced from \citer\hasse(Kap.~16) or
\citer\course(Lecture~18), and equals the number of cyclic subgroups of
$F^\times\!/F^{\times(p-1)}$.  For a given $F'$, a cyclic degree-$p$ extension
$L|F'$ arises from some $E|F$ if and only if $L|F$ is galoisian, and $E_1,E_2$
give rise to the same $L$ if and only if they are conjugate over $F$
(cf.~lemma.~20).

But it is perhaps simpler to work over the compositum
$K=F(\!\root{p-1}\of{F^\times})$ of all $F'$~; one advantage is that $\root
p\of1\in K$ (true by convention in characteristic~$p$).

{\bf 3. The compositum of all separable extensions of degree~$p$}\pointir We
are in a position to the main theorem of \citer\delcorso(th.~1) prove the
characteristic-$p$ analogue.

\th COROLLARY 19
\enonce
Every degree-$p$ separable extension of\/ $F$ becomes cyclic when translated
to\/ $K$.
\endth 
Indeed, as $F$ contains the $(p-1)$-th roots of~$1$, the compositum of all
abelian extensions of exponent dividing $p-1$ is $K$.

\smallskip

Put $G=\Gal(K|F)$~; it is a free $(\Z/(p-1)\Z)$-module of rank~2, dual to
$F^\times\!/F^{p-1}$.

Which cyclic extensions $L|K$ arise as $L=EK$ for some (degree-$p$, separable)
extension $E|F$~?  If $L$ does, then $L|F$ would be galoisian, for there is
some cyclic extension $F'$ of $F$ in $K$ such that $EF'|F$ is galoisian
(prop.~??).  Conversely,

\th LEMMA 20
\enonce
If\/ $L$ is galoisian over $F$, then there is a degree-$p$ separable
extension\/ $E|F$ such that\/ $L=EK$~; two such extensions $E,E'$ give rise to
the same $L$ if and only if they are conjugate over\/ $F$. 
\endth
Let $H=\Gal(L|K)$, and suppose that $L|F$ is galoisian, so that there is an
exact sequence $1\to H\to\Gal(L|F)\to G\to1$.  The action of $G$ on $H$ by
cojugation corresponds to some character $\chi:G\to\F_p^\times$, so that for
every $\sigma\in H$ and every $\tau\in\Gal(L|F)$, we have
$\tau\sigma\tau^{-1}=\sigma^{\chi(\bar\tau)}$, where $\bar\tau\in G$ is the
image of $\tau$.  We also have $\Gal(L|F)=H\times_\chi G$, which has at least
one index-$p$ subgroup, for example $\{1\}\times_\chi G$, so there is at least
one degree-$p$ separable extension $E|F$ such that $L=EK$.  We have already
seen (prop.~18) that any two such $E$ are conjugate over $F$.

\smallskip

But the great things about degree-$p$ cyclic extensions $L|K$ is that they
correspond to lines in the $\Fp$-space $V=K^\times\!/K^{\times p}$ in the case
$F|\Qp$, and in $V=K/\wp(K)$ in the case $F|\Fp(\!(\pi)\!)$.  When is the
(degree-$p$, cyclic) extension corresponding to a line $D\subset V$ galoisian
over $F$~?  Precisely when $D$ is stable under the action of $G$ on $V$
(cf.~lemmas~1 and~2)~:

\th PROPOSITION 21
\enonce
If the line\/ $D\subset V$ and the degree-$p$ cyclic extension\/ $L|K$
correspond to each other, then\/ $L|F$ is galoisian if and only if\/ $D$ is\/ 
$G$-stable.  If so, $L|K$ is abelian (and indeed cyclic) if and only if the
action of\/ $G$ on\/ $D$ is via\/ $\omega$.
\endth

Notice that if $G$ acts on $D$ via $\chi:G\to\F_p^\times$, then the
conjugation action of $G$ on (the order-$p$ group) $\Gal(L|K)$ is via
$\chi\omega^{-1}$, so $\Gal(L|F)$ is commutative if and only if
$\chi=\omega$. 

Suppose that $D\subset V$ is a $G$-stable line and that $G$ acts on it by
$\chi:G\to\F_p^\times$.  Put $L=K(\root p\of D)$ (resp.~$LK(\wp^{-1}(D)$, so
that $L|F$ is galoisian~; put $H=\Gal(L|K)$.  we have an exact sequence
$$
1\to H\to\Gal(L|F)\to G\to1.
$$
which gives an action of $G$ on $H$ by conjugation.

\th PROPOSITION
\enonce
If the action of\/ $G$ on\/ $D$ is via the character\/ $\chi:G\to\F_p^\times$,
then its action on\/ $H$ is via the character\/ $\chi\omega^{-1}$.
\endth

Consider the extension $K(\root p\of{K^\times})$ in the characteristic-$0$
case, and the extension $K(\wp^{-1}(K))$ in the characteristic-$p$ case.

\th COROLLARY 22
\enonce
The extension\/ $K(\root p\of{K^\times})$ (resp.~$K(\wp^{-1}(K))$) is the
compositum of all degree-$p$ separable extensions of\/ $F$. 
\endth
This is the main result of \citer\delcorso() in the characteristic-$0$ case.
As $K(\root p\of{K^\times})$ (resp.~$K(\wp^{-1}(K))$) is the maximal abelian
extension of exponent~$p$, it is sufficient, in view of the preceding
statements, to show that $V$ is generated by the union of the $G$-stable lines
it contains.  This is true because the $\Fp[G]$-module $V$ is the internal
direct sum, over the various characters $\chi:G\to\F_p^\times$, of the
$\chi$-eignespaces $V(\chi)\subset V$ for the action of $G$ on $V$, thanks to
the fact that the field $\Fp$ contains all $(p-1)$-th roots of~1 (Fermat) and
that the group $G$ is commutaive of exponent $p-1$.

{\bf 4.  Filtered galoisian modules}\pointir 

Examples of $G$-stable lines are provided, in the case of $F|\Qp$ of
ramification index $e$, by $\bar U_{pe}$ and $\bar\mu$ (the image of the
torsion subgroup $\mu\subset{\goth o}^\times$), on both of which $G$ acts
via the cyclotomic character $\omega$.  In the characteristic-$p$ case, the
line $\bar{\goth o}=\Fp$ is stable and the action of $G$ is in fact trivial.

The subspaces $\bar U_i$ for $i\in[0,pe]$ (resp.~$\overline{{\goth
    p}^i}$ for $i\in-\N$) are $G$-stable, essentially because there is a
unique extension of the valuation from $F$ to $K$.  We have seen that $\bar
U_{pi+1}=\bar U_{pi}$ except for $i=0,e$ (\citer\locdisc(prop.~42)) and that
$\overline{\pgoth^{pi+1}}=\overline{\pgoth^{pi}}$ except for $i=0$
(\citer\further(prop.~11)).

The codimension is~$1$ in the three exceptional cases.  In characteristic~$0$,
we have $\bar U_{pe+1}=\{1\}$ and $\bar U_{pe}$ is a stable line on which $G$
acts via $\omega$, and the valuation $v_K$ provides an isomorphism $\bar
U_0/\bar U_1\to\Z/p\Z$.  In characteristic~$p$, we have $\bar\pgoth=0$ and the
trace $S_{l|\Fp}$ provides an isomorphism $\bar{\goth o}\to\Fp$, where $l$
is the residue field of $K$.  For all other $i$, the codimension equals the
absolute degree $[l:\Fp]$, and indeed the quotients are canonically
isomorphic to $U_i/U_{i+1}$ (resp.~${\goth p}^{i}/{\goth p}^{i+1}$)
(\citer\locdisc(prop.~42), (\citer\further(prop.~11)).  Thus they are not
merely $\Fp$-spaces but $k$-spaces.  The pictures in \citer\final(), \S5
summarise some of these facts.

\th PROPOSITION 23
\enonce
For these\/ $i$, the isomorphisms\/ $\bar U_i/\bar U_{i+1}\to U_i/U_{i+1}$
(resp.~$\overline{\pgoth^i}/\overline{\pgoth^{i+1}} \to
\pgoth^i/\pgoth^{i+1}$) of\/ $k$-spaces are $G$-equivariant.  In the
former case, the $k$-isomorphisms\/ $U_i/U_{i+1}\to \pgoth^i/\pgoth^{i+1}$ are 
also\/ $G$-equivariant. 
\endth
Indeed, the isomorphisms in the first assertion are induced by the identity
map on $U_i$ (resp.~${\goth p}^i$)~: the diagrams
$$
\def\\{\mskip-2\thickmuskip}
\def\droite#1{\\\hfl{#1}{}{8mm}\\}
\def\vide{\phantom{phantom}}
\diagram{
U_i&\droite{=}&U_i&\vide&\pgoth^i&\droite{=}&\pgoth^i\cr
\vfl{}{}{5mm}&&\vfl{}{}{5mm}&\vide&\vfl{}{}{5mm}&&\vfl{}{}{5mm}\cr
\bar U_i/\bar U_{i+1}&\droite{\sim}&U_i/U_{i+1}
 &\vide&\overline{\pgoth^i}/\overline{\pgoth^{i+1}}&\droite{\sim}
 &\pgoth^i/\pgoth^{i+1}\cr
}
$$ 
are commutative.  The second assertion is also clear, as the isomorphism in
question is nothing but $\bar u\mapsto\overline{u-1}$.

\smallskip

In order to study the $G$-modules $\pgoth^i/\pgoth^{i+1}$, a preliminary
study of the $G$-module $l$ is fruitful.

\th LEMMA ??
\enonce
Let\/ $k$ be a finite field, $q=\Card k$, and\/ $l|k$ any extension of degree 
dividing $q-1$.  For every character\/ $\chi:\Gal(l|k)\to k^\times$, the
$\chi$-eigenspace $l(\chi)$ is a $k$-line in $l$.
\endth
Let $\varphi$ (Frobenius) be the canonical generator $x\mapsto x^q$ of
$\Gal(l|k)$, and put $a=\chi(\varphi)$~; let $m=[l:k]$ be the order of
$\varphi$, so that the order of $a\in k^\times$ divides~$m$.  The
$\chi$-eigenspace consists of all $x\in l$ such that $\varphi(x)=ax$~; such
$x$ are roots of the binomial $T^q-aT=T(T^{q-1}-a)$, which has at most $q$
roots.  As it defines a linear endomorphism of the $k$-space $l$, it is
sufficient to prove that $a$ has a $(q-1)$-th root in $l$.

We have said that the order of $a$ in $k^\times\!/k^{\times q-1}=k^\times$
divides $m$, Therefore the degree of the extension $k(\!\root{q-1}\of a)|k$
divides $m$, and hence $k(\!\root{q-1}\of a)\subset l$.  This shows that the
$k$-endomorphism $x\mapsto x^q-ax$ of $l$ is not injective, and hence its
kernel $l(\chi)$ is a $k$-line in $l$.  Incidentally, if $\chi$ is trivial,
then $a=1$ and $l(\chi)=k$.
 
\smallskip

Momentarily let $K$ be any galoisian extension of $F$ of group $G$, and
suppose that there is a uniformiser $\pi$ of $K$ such that $\varpi=\pi^s$ is
in $F$ for some $s>0$.

\th PROPOSITION 25
\enonce
For every integer\/ $i\in\Z$, ``\thinspace multiplication by\/
$\varpi$'' gives an isomorphism\/
$\pgoth^i\!/\pgoth^{i+1}\to\pgoth^{i+s}\!/\pgoth^{i+s+1}$ of\/ $k[G]$-modules.
\endth
More precisely, the reduction modulo $\pgoth$ of the $\ogoth$-linear
isomorphism $x\mapsto\varpi x:\pgoth^i\to\pgoth^{i+s}$ is $G$-equivariant.
But this is clearly the case~: $\sigma(\varpi x)=\varpi\sigma(x)$ for every 
$\sigma\in G$, because $\varpi\in F$ and $\sigma$ is $F$-linear.

\smallskip

Let us now return to our $K=F(\!\root{p-1}\of{F^\times})$ and
$G=\Hom(F^\times\!/F^{\times p-1},\F_p^\times)$, so that the group of
characters of $G$ is $\Hom(G,\F_p^\times)=F^\times\!/F^{\times p-1}$.  (The
cyclotomic character $\omega$ corresponds to $\overline{-p}$ in
characteristic~$0$, and $\omega=1$ by convention in characteristic~$p$).  Each
character $\chi$ therefore has a ``\thinspace valuation\thinspace'' $\bar
v(\chi)\in\Z/(p-1)\Z$, coming from the valuation $v:F^\times\to\Z$.
Unramified characters --- those in the kernel
$\ogoth_F^\times\!/\ogoth_F^{\times p-1}$ of $\bar v$ --- are the same as
characters of the quotient $G/G_0=\Gal(l|k)$, where $G_0\subset G$ is the
inertia subgroup.  In characteristic~$0$, the cyclotomic character $\omega$ is
unramified if and only if $p-1|e$, for $\omega$ corresponds to $\overline{-p}$
and $\bar v(\overline{-p})\equiv e\pmod{p-1}$.

Let us decompose the $k[G]$-modules $\pgoth^i/\pgoth^{i+1}$ as an internal
direct sum of $\chi$-eigenspaces for various $\chi:G\to\F_p^\times$.  For
$i=0$, the $k[G]$-module $\pgoth^i/\pgoth^{i+1}$ is in fact the
$k[G/G_0]$-module $l$.  We have seen that for every unramified character
$\chi$ of $G$, the $\chi$-eigenspace is a $k$-line (lemma~24).  It follows
that for every {\it ramified\/} character $\chi$, we have $l(\chi)=0$~: there
is room only for so many, and unramified characters have used it all up.

Let us provide the details of the notion of twisting a $k[G]$-module by a
character $\xi:G\to\F_p^\times$.  Denote by $l\{\xi\}$ the $k[G]$-module
whose underlying $k$-space is $l$, but the new action $\star_\xi$ is defined
by $\sigma\star_\xi x=\xi(\sigma)\sigma(x)$ for every $\sigma\in G$ and every
$x\in l$, so that if $\xi=1$ is the trivial character, then $l\{\xi\}$ is
the $k[G]$-module $\ogoth/\pgoth$.  In this process, the $\chi$-eigenspaces of
$l$ get converted into $\xi\chi$-eigenspaces in $k\{\xi\}$.

It is easy to see that $l\{\xi\}$ depends only on $\bar
v(\xi)\in\Z/(p-1)\Z$ (up to $k[G]$-isomorphism).  We denote by $l[i]$ (for
$i\in\Z$) the $k[G]$-modules $l[\xi]$ for any $\xi$ such that $\bar
v(\xi)\equiv i\pmod{p-1}$.

\th PROPOSITION 26
\enonce
For every\/ $i\in\Z$ and every character\/ $\chi:G\to\F_p^\times$, the\/
$\chi$-eigenspace in the\/ $k[G]$-module\/ $l[i]$ is a\/ $k$-line if\/ $\bar
v(\chi)\equiv i\,(p-1)$~; it is reduced to\/~$0$ otherwise.
\endth
This is just lemma~24 in the case~$i=0$, and the general case follows from
this by our discussion of twisting.  An immediate consequence is the following
result. 

\th PROPOSITION 27
\enonce
For\/ $i\in\Z$, the\/ $k[G]$-module\/ $\pgoth^i\!/\pgoth^{i+1}$ is
isomorphic to\/ $l[i]$.
\endth
Choose any uniformiser $\pi$ of $K$ such that $\pi^{p-1}$ is (a uniformiser)
in $F$~; this is possible.  It is easy to see that, by taking $\pi^i$ as an
$\ogoth$-basis of $\pgoth$, the resulting $k$-linear map
$\pgoth^i/\pgoth^{i+1}\to l\{\xi^i\}$ is an isomorphism of $k[G]$-modules,
where $\xi$ is the character such that $\sigma(\pi)=\xi(\sigma)\pi$ for every
$\sigma\in G$.  As we have $\bar v(\xi)\equiv1$, this
shows that $\pgoth^i/\pgoth^{i+1}$ is $k[G]$-isomorphic to $l[i]$.

\th COROLLARY 28
\enonce
Put\/
$W=\pgoth^{-1}\!/\pgoth^{-2}\oplus\cdots 
\oplus\pgoth^{-(p-1)}\!/\pgoth^{-p}$. 
For every\/ $\chi:G\to\F_p^\times$, the\/ $\chi$-eigenspace in
the\/ $k[G]$-module\/ $W$ is a\/ $k$-line. 
\endth
In other words, we have endowed $W$ with the filtration for which the
successive quotients are, for example when $p=5$, 
$$
l[-1],\ l[-2],\ l[-3],\ l[-4],
$$
in that specific order, rather than in some other order.  If we twist it by
the cyclotomic character (which we will soon need to) to get $W\{\omega\}$,
and if $\omega$ has ``\thinspace valuation\thinspace'' $\bar v(\omega)=2$,
then the successive quotients of the filtered $\Fp[G]$-module $W\{\omega\}$
are
$$
l[1],\ l[0],\ l[-1],\ l[-2],
$$
which are shifted two steps to the right (or to the left, as
$2\equiv-2\pmod4$).  We need to keep track of both the filtration and the
$G$-action.  There is no difference when $\omega$ is unramified, for then the
shift is by~$0$ steps. 

Get back to the $k[G]$-modules $\bar U_i/\bar U_{i+1}$
(resp.~$\overline{\pgoth^i}/\overline{\pgoth^{i+1}}$) for appropriate~$i$.

\th PROPOSITION 29
\enonce
The\/ $k[G]$-module\/ $\bar U_i/\bar U_{i+1}$ for\/ $0<i<pe$ prime to\/~$p$ in
the characteristic-$0$ case
(resp.~$\overline{\pgoth^i}/\overline{\pgoth^{i+1}}$ for\/ $i<0$ prime
to\/~$p$ in the characteristic-$p$ case) is isomorphic to\/  $l[i]$.
\endth

Let us record all this in a picture, for I have still not got over the fact
that the quotients are $k[G]$-modules for appropriate~$i$, instead of merely
being $\Fp[G$]-modules.  In characteristic~$0$, the picture begins with
$\displaystyle
\{\bar1\}\underbrace{\subset}_{\Fp\{\omega\}}\bar U_{pe}
\underbrace{\subset}_{l[pe-1]}\bar U_{pe-1}\cdots
$
to end with
$$
\cdots
\underbrace{\subset}_{l[pi+1]}\bar U_{pi+1}
=\bar U_{pi}
\underbrace{\subset}_{l[pi-1]}\bar U_{pi-1}\cdots
\underbrace{\subset}_{l[1]}\bar U_{1}
\underbrace{\subset}_{\Fp}\bar U_{0},
$$
whereas in characteristic~$p$ it goes on for ever
$$
\{\bar0\}\underbrace{\subset}_{\Fp}
\bar\ogoth\underbrace{\subset}_{l[-1]}
\overline{\pgoth^{-1}}
\cdots
\underbrace{\subset}_{l[pi+1]}\overline{\pgoth^{pi+1}}
=\overline{\pgoth^{pi}}
\underbrace{\subset}_{l[pi-1]}
\overline{\pgoth^{pi-1}}
\cdots
$$
The beauty of this can reduce even the most hardened criminal to tears.

The analogy can be further improved.  First, declare $\omega$ to be the
trivial character in characteristic-$p$ case.  Second, the two pictures will
look even more similar if the first one is shifted to the right by $pe$.  The
problem is that the $k[G]$-modules $l[pe-1]$ and $l[-1]$ are not isomorphic,
unless $\omega$ is unramified.  More precisely, if $s\in[0,p-1[$ is the
integer such that $e\equiv s\pmod{p-1}$, so that $\bar v(\omega)=s$, then
$l[pe-1]$ is the sum of $k$-lines indexed by the characters of ``valuation''
$s-1$, whereas $l[-1]$ is the sum of $k$-lines indexed by the characters of
``valuation'' $-1$.  But this can be easily remedied if we twist the latter by
$\omega$.  

So $l[-1]\{\omega\}$ is the same $k[G]$-module as $l[pe-1]$.  Suppressing
the terms indexed by multiples of~$p$, and exploiting the fact that any
$G$-stable subspace has a unique $G$-stable supplement, the $\Fp[G]$ module
$K^\times\!/K^{\times p}$ (resp.~$K/\wp(K)$) is
$$
\Fp\{\omega\}\oplus l[-1]\{\omega\}\oplus\cdots 
\oplus l[-b^{(j)}]\{\omega\}\oplus\cdots\ \big(\oplus\Fp\big)
$$
where the middle terms are indexed by the sequence $b^{(j)}$ of
\hbox{prime-to-$p$} integers, for every $j>0$ in characteristic~$p$ but only
for $j\in[1,(p-1)e]$ in characteristic~0, which is also when the last term
appears.  This picture keeps track of both the filtration and the $G$-action. 

Grouping together $p-1$ middle terms at a time, and recalling that
$W=\pgoth^{-1}\!/\pgoth^{-2}\oplus\cdots \oplus\pgoth^{-(p-1)}\!/\pgoth^{-p}$
by definition (cor.~28), we get our final picture of the $\Fp[G]$-module
$K^\times\!/K^{\times p}$ (resp.~$K/\wp(K)$), namely
$$
\Fp\{\omega\}\oplus W\{\omega\}\oplus W\{\omega\}\oplus\cdots\
\big(\oplus\Fp\big)~;
$$
there are $e$ copies of $W\{\omega\}$ in the characteristic-$0$ case and
the final parenthetical term appears, but infinitely many copies of
$W\{\omega\}$ in the characteristic-$p$ case (and the final term does not
appear).  The twist by $\omega$ doesn't change anything if $\omega$ is
unramified.  This picture carries the same information as the previous ones,
only less visibly 
so.  

{\bf 5. Serre's mass formula in prime degree}\pointir 

It is time to do the counting.  Let $E|F$ be a ramified separable extension of
degree~$p$ and let $D\subset K^\times/K^{\times p}$ (resp.~$D\subset
K/\wp(K)$) be the $G$-stable line corresponding to $E$.  Using the {\it
  Schachtelungssatz}, it is easy to see that $c(E)=b$, where $b$ is the unique
ramification break of $\Gal(K(\root p\of D)|K)$.  Also, the number of $E$
which give rise to the same $D$ is~$1$ if the character $\chi$ through which
$G$ acts on $D$ is $\omega$, and $p$ if $\chi\neq\omega$.  In short, the
contribution of such $E$ to Serre's mass formula is
$$
\sum_{E\mapsto D}q^{-c(E)}=\cases{
\phantom{p}q^{-b}&if $\chi=\omega$\phantom.\cr
pq^{-b}&if $\chi\neq\omega$.\cr}
$$
So the sum over all $E$ gets replaced by a sum over all~$D$. 

The dimension of the space $\chi$-eigenspace $(\Fp\{\omega\}\oplus
W\{\omega\}^i)(\chi)$ (for all $i\in\N$ in characteristic~$p$, for $i\in[0,e]$
in characteristic~$0$) is
$$
\dim_{\Fp}(\Fp\{\omega\}\oplus W\{\omega\}^i)(\chi)
=\cases{
1+if&if $\chi=\omega$\phantom.\cr
\phantom{1+\ }if&if $\chi\neq\omega$.\cr}
$$
So the number of lines in $(\Fp\{\omega\}\oplus W\{\omega\}^{i+1})(\chi)$
which are not in $(\Fp\{\omega\}\oplus W\{\omega\}^i)(\chi)$ is
$$
{pq^{i+1}-1\over p-1}-{pq^{i}-1\over p-1}\qquad
\left(\hbox{resp. }{q^{i+1}-1\over p-1}-{q^{i}-1\over p-1}\right).
$$
Assume for a moment that $\omega=1$, as is always the case in
characteristic~$p$.  The unique ramification break of such a line depends only
on $i$ and the ``\thinspace valuation\thinspace'' $j\in[1,(p-1)]$ of $\chi$
and equals $pi+j$.  The contribution of all such lines (for a
given~$i$ and $\chi$ of valuation~$j$) is 
$$
\left({pq^{i+1}-1\over p-1}-{pq^{i}-1\over p-1}\right)q^{-(pi+j)}
$$
irrespective of whether $\chi=\omega$ or $\chi\neq\omega$.  As there are $p-1$
characters for a given~$j$, the contribution of all $G$-stable lines for a
given $i$ is 
$$
p(q^{i+1}-q^{i})q^{-pi}Q,\qquad Q=\sum_{j=1}^{p-1}q^{-j}.
$$
Now all that remains to do in the characteristic-$p$ case is to sum over all
$i\in\N$ and use the fact that
$$
\sum_i(q^{i+1}-q^i)q^{-pi}=(q-1)\sum_i q^{i-pi}=(q-1){q^{p-1}\over
  q^{p-1}-1}={1\over Q},
$$
proving the formula in the characteristic-$p$ case.  In the characteristic-$0$
case, the sum extends only over $i\in[0,e]$ but the contribution of the {\it
  tr{\`e}s ramifi{\'e}es\/} extensions, which correspond to (all) lines in
$\Fp\{\omega\}\oplus W\{\omega\}^e\oplus\Fp$ which are not in
$\Fp\{\omega\}\oplus W\{\omega\}^e$, should be added, for they are in the
$1$-eigenspace.  The reader will have no difficulty in verifying that the sum
turns out to be $p$ in this case as well.

We had assumed that $\omega$ is trivial.  What if $\omega\neq1$~?  It is clear
that the contribution of $\omega$ is the same as the amount by which the
contribution of the trivial character goes down~: the sum of the two
contributions is the same.  Also, the sum of the contributions from the other
characters remains the same~: for each~$i$, the $p-1$ terms, each coming from
to the $p-1$ characters with a given ``\thinspace valuation\thinspace''
$j\in[1,p-1]$, get cyclically permuted by the addition of $\bar v(\omega)$.
It all balances out.  (The contribution of $\omega$ can be independently
computed, see \citer\final(),~\S6.)

We have proved the degree-$p$ case of Serre's mass
formula~:

\th THEOREM 30 (Serre, 1978)
\enonce
Let\/ $F$ be a local field with finite residue field of characteristic~$p$
and cardinality~$q$.  When\/ $E$ runs through ramified separable degree-$p$
extensions of\/ $F$,  
$$
\sum_E q^{-c(E)}
=p,
$$
where\/ $c(E)=v(d_{E|F})-(p-1)$ and $v(d_{E|F})>(p-1)$ is the valuation of
the  discriminant\/ $d_{E|F}$ of\/ $E|F$.
\endth

\bigbreak
\unvbox\bibbox 

\bye